\documentstyle[11pt,twoside,amsfonts,amssymb]{article}
  
\begin{titlepage}

\title{$L^2$-invariants of locally symmetric spaces}

\author{Martin
Olbrich\thanks{Mathematisches Institut, Universit\"at G\"ottingen, Bunsenstr. 3-5, 37073 G\"ottingen, GERMANY, E-mail: olbrich@uni-math.gwdg.de}
}

\end{titlepage}

\newcommand{\proof}{{\it Proof.$\:\:\:\:$}}
\newcommand{\D}{\displaystyle}

\newcommand{\kaaa}{{\frak k}}
\newcommand{\paaa}{{\frak p}}
\newcommand{\taaa}{{\frak t}}
\newcommand{\haaa}{{\frak h}}

\newcommand{\R}{{\Bbb R}}

\newcommand{\C}{{\Bbb C}}

\newcommand{\gaaa}{{\frak g}}
\newcommand{\maaa}{{\frak m}}
\newcommand{\aaaa}{{\frak a}}
\newcommand{\naaa}{{\frak n}}

\newcommand{\cZ}{{\cal Z}}

\newcommand{\cC}{{\cal C}}

\newcommand{\vol}{{\mathrm{ vol}}}

\newcommand{\ve}{\varepsilon}

\newcommand{\End}{{ \mathrm{ End}}}

\newcommand{\tr}{{ \mathrm{ tr}}}
\newcommand{\Tr}{{ \mathrm{ Tr}}}

\newcommand{\Ad}{{\mathrm{ Ad}}}

\newcommand{\id}{{ \mathrm{ id}}}

\newcommand{\aca}{{\aaaa_\C^\ast}}

\def\hB{\hspace*{\fill}$\Box$\newline\noindent}

\newtheorem{prop}{Proposition}[section]
\newtheorem{lem}[prop]{Lemma}

\newtheorem{theorem}[prop]{Theorem}
\newtheorem{kor}[prop]{Corollary}

\def\imath{i}

\begin{document}
\setcounter{page}{1}
\maketitle

\begin{abstract}
Let $X=G/K$ be a Riemannian symmetric space of the noncompact type,
$\Gamma\subset G$ a discrete, torsion-free, cocompact subgroup, and let
$Y=\Gamma\backslash X$ be the corresponding locally symmetric space.
In this paper we explain how the Harish-Chandra Plancherel Theorem for
$L^2(G)$ and results on
$({\frak g}, K)$-cohomology can be used in order to compute the $L^2$-Betti numbers, the Novikov-Shubin invariants, and the
$L^2$-torsion of $Y$ in a uniform way thus completing results previously
obtained by Borel, Lott, Mathai, Hess and Schick. It turns out that the behaviour of these invariants
is essentially determined by the fundamental
rank $m=\mbox{rk}_{\C}G- \mbox{rk}_{\C}K$ of $G$. In particular, we show
the nonvanishing of the $L^2$-torsion of $Y$ whenever $m=1$.
\end{abstract}

\tableofcontents
 
\parskip3ex

\section{Introduction}

During the last two decades $L^2$-invariants have proved to be a powerful tool in the topology
of compact manifolds (see \cite{lueck97} for an overview). Although they can be defined in purely combinatorial terms we are interested here in their equivalent analytic versions: They are spectral invariants of the $p$-form Laplacians of
the universal cover of the manifold. For particular nice manifolds these might be
computable. Indeed, the aim of the present paper is to extract from the representation theoretic work of Harish-Chandra
\cite{HC76III} and Borel-Wallach \cite{borelwallach80} information on the spectral decomposition of the form Laplacians on Riemannian symmetric
spaces of the non-compact type which is sufficiently explicit in order to compute
the spectral invariants of interest. We try to do this in a rather
detailed way which, we hope, keeps the paper readable for nonspecialists
in harmonic analysis.

Let $X\rightarrow Y$ be the universal cover of a compact Riemannian manifold. Set $\Gamma:=\pi_1(Y)$. The form Laplacian $\Delta_p=d^*d+dd^*$ defines a non-negative, elliptic,
self-adjoint operator acting on $L^2(X,\Lambda^pT^*X)$, the square integrable $p$-forms on $X$. By
$\Delta^\prime_p$ and $\Delta^c_p=(d^*d)^\prime$ we denote the restriction of $\Delta_p$ to the orthogonal
complement of its kernel and the coclosed forms in this orthogonal complement, respectively. We consider the corresponding heat kernels $e^{-t\Delta_p^*}(x,x^\prime):=(P_*e^{-t\Delta_p})(x,x^\prime)$, $x,x^\prime \in X$, for $*=\emptyset,\prime$ or $c$. Here $P_*$ denotes the orthogonal projection to the corresponding subspace. The local traces
$\tr\, e^{-t\Delta_p^*}(x,x)$ are $\Gamma$-invariant functions on $X$. For a thorough discussion
of the following definitions we refer to \cite{lott92}, \cite{mathai92}, and \cite{lueck97}.

We set
$$ \Tr_\Gamma e^{-t\Delta_p^*}:=\int_F \tr\, e^{-t\Delta_p^*}(x,x)\:dx\ ,$$
where $F\subset X$ is a fundamental domain of the action of $\Gamma$ on $X$ and $dx$
is the Riemannian volume element of $X$. Then the $L^2$-Betti numbers are given
by
$$ b_p^{(2)}(Y):=\lim_{t\to\infty} \Tr_\Gamma e^{-t\Delta_p}\in [0,\infty)\ .$$
They are equal to the von Neumann dimension of $\ker \Delta_p$ viewed as
Hilbert ${\cal N}(\Gamma)$-module, where ${\cal N}(\Gamma)$ is the group von Neumann
algebra of $\Gamma$. If the spectrum of $\Delta_p$ has no gap around $0$
the Novikov-Shubin invariants of $Y$ are defined by
$$ \tilde \alpha_p(Y):=\sup\{\beta\:|\: \Tr_\Gamma e^{-t\Delta_p^\prime}=O(t^{-\frac{\beta}{2}})\mbox{ as }t\to\infty\}\in [0,\infty]\ .$$
It measures the asymptotic behaviour of the spectral density function of $\Delta_p$ at $0$. In case of a gap around $0$ we set $\tilde\alpha_p(Y):=\infty^+$. Replacing $\Delta^\prime_p$ by $\Delta^c_p$ we
obtain the analogously defined Novikov-Shubin invariants $\alpha_p(Y)$ of
$d^*d$. Using the action of the exterior differential $d$ on the Hodge decomposition of $L^2$-forms we obtain \begin{equation}\label{lott}
\tilde\alpha_p(Y)=\min\{\alpha_p(Y),\alpha_{p-1}(Y)\}\ .
\end{equation}
Finally, if $\alpha_p(Y)>0$ for all $p$ (or, more generally, if $X$ is of determinant class (see \cite{lueck97})), then the $L^2$-torsion of $Y$ is defined by
$$ \rho^{(2)}(Y):=\frac{1}{2}\sum (-1)^{p+1}p\:\log{\det}_\Gamma(\Delta^\prime_p)=\frac{1}{2}\sum (-1)^{p}\log{\det}_\Gamma(\Delta^c_p)\ ,$$
where for $*=\prime, c$
$$ -\log{\det}_\Gamma(\Delta_p^*):=\frac{d}{ds}_{|s=0}\left(\frac{1}{\Gamma(s)}\int_0^\ve
\Tr_\Gamma e^{-t\Delta_p^*}t^{s-1}dt\right)+\int_\ve^\infty \Tr_\Gamma e^{-t\Delta_p^*}t^{-1}dt $$
for any $\ve>0$, where the first integral is considered as a meromorphic function in $s$. By Poincar\'e duality $\rho^{(2)}(Y)=0$ for even dimensional
manifolds $Y$.

From now on let $X=G/K$ be a Riemannian symmetric space of the noncompact type.
Here $G$ is a real, connected, linear, semisimple Lie group without compact
factors, and $K\subset G$ is a maximal compact subgroup. It is the universal
cover of compact locally symmetric spaces of the form $Y=\Gamma\backslash X$,
where $\Gamma\cong \pi_1(Y)$ can be identified with a discrete, torsion-free,
cocompact subgroup of $G$. It will be convenient to consider also the compact
dual $X^d$ of $X$. $X^d$ is defined as follows: Let
$\gaaa$, $\kaaa$ be the Lie algebras of $G$, $K$. Then we have the Cartan
decomposition $\gaaa=\kaaa\oplus\paaa$. Then $\gaaa^d:=\kaaa\oplus i\paaa$ is
another subalgebra of the complexification of $\gaaa$. Let $G^d$ be the corresponding analytic subgroup of the complexification $G_\C$ of $G$. Then
$G^d$ is a compact group, and $X^d=G^d/K$. We normalize the Riemannian
metric on $X^d$ such that multiplication by $i$ becomes an isometry $T_{eK}X\cong\paaa\rightarrow i\paaa\cong T_{eK}X^d$. In this paper we are going to prove the following theorem:

\begin{theorem}\label{main} 
Let $n=\dim Y$ and $m=m(X):={\mathrm rk}_{\Bbb C}G- {\mathrm rk}_{\Bbb C}K$ be the fundamental
rank of $G$. Let $\chi(Y)$ be the Euler characteristic of $Y$. Then
\begin{enumerate}
\item[(a)] $b^{(2)}_p(Y)\ne 0\Leftrightarrow m=0$ and $p=\frac{n}{2}$. In particular, $b^{(2)}_\frac{n}{2}(Y)=(-1)^\frac{n}{2}\chi(Y)=\D\frac{\vol(Y)}{\vol(X^d)}\chi(X^d)$.
\item[(b)] $\alpha_p(Y)\ne \infty ^+\Leftrightarrow m>0$ and $p\in [\frac{n-m}{2},\frac{n+m}{2}-1]$.
In this range $\alpha_p(Y)=m$.
\item[(c)] $\rho^{(2)}(Y)\ne 0\Leftrightarrow m=1$.
\end{enumerate}
\end{theorem} 
Note that $n-m$ is always even and positive.
Part (a) of the theorem was known for a long time, at least since Borel's
paper \cite{borel85}. For special cases see also \cite{dodziuk79}, \cite{donnelly81}. For the convenience of the reader we include a proof here. 
In fact, we prove a stronger statement which should have been known to
the experts although we were not able to find it in the literature:

\begin{prop}\label{spec}
The discrete spectrum of $\Delta_p$ on $L^2(X,\Lambda^pT^*X)$ is empty unless $m(X)=0$ and $p=\frac{n}{2}$.
In this case $0$ is the only eigenvalue of $\Delta_p$.
\end{prop}

The strategy of the proof of (b) can already be found in Lott's paper \cite{lott92}, Section VII. To be more precise, Equation (\ref{lott}) 
implies the slightly weaker result
\begin{equation}\label{lm}
\tilde \alpha_p(Y)=\left\{\begin{array}{cl}
\infty^+& p\not\in [\frac{n-m}{2},\frac{n+m}{2}] \mbox{ or } m=0\\
m       & p\in [\frac{n-m}{2},\frac{n+m}{2}] \mbox{ and } m\ne 0
\end{array}\right. \ . 
\end{equation}
Lott proved the first line of (\ref{lm}) and that $\tilde \alpha_p(Y)$ is finite
and independent of $p$ for the remaining values of $p$. He indicated how one should be able
to compute the precise value of $\tilde \alpha_p(Y)$. But he finished the computation in the real rank one case, only. In addition, already Borel \cite{borel85} showed
that the range of the differential $d_p$ of the $L^2$-de-Rham-complex is 
not closed
for $p\in [\frac{n-m}{2},\frac{n+m}{2}-1]$. Theorem \ref{main} (b) can be
considered as a quantitative refinement of this result. After the present paper
was written I was informed by S. Mehdi that there is a joint paper of him
with  N. Lohoue \cite{lohouemehdi00} which has recently appeared in print and
which contains a proof of (\ref{lm}). The interested reader will also find 
there more information concerning the material presented here in Sections \ref{shc} and \ref{coho}. But he should be aware that in that paper the range 
of finiteness
of $\tilde \alpha_p(Y)$ is constantly misprinted and that the proof of
the second line of (\ref{lm}) as written down there is not quite complete (it is
not mentioned that it is important to know that $p_\xi(0)>0$, see Equation
(\ref{crux}) below).
    
The main motivation to do the present work was to obtain part (c) of the theorem.
In our locally homogeneous situation we have for any $x\in X$
$$ \Tr_\Gamma e^{-t\Delta_p^*}=\vol(Y)\cdot \tr\, e^{-t\Delta_p^*}(x,x)$$
and thus
$$\rho^{(2)}(Y)=\vol(Y)\cdot T^{(2)}(X)$$
for a certain real number $T^{(2)}(X)$. Note that in contrast to 
$\rho^{(2)}(Y)$ the number $T^{(2)}(X)$ depends on the normalization of the invariant metric on $X$ (of course only via the volume form).
A well-known symmetry argument (\cite{moscovicistanton91}, Proposition 2.1) yields that $T^{(2)}(X)=0$ whenever $m\ne 1$. Lott \cite{lott92} (see also \cite{mathai92})
showed that $T^{(2)}(H^n)\ne 0$ for $n=3,5,7$, where $H^n$ is the real hyperbolic
space (his values for $T^{(2)}(H^n)$ for $n=5,7$ were not correct).   
This led to the conjecture that $T^{(2)}(H^n)\ne 0$ for all odd $n$ which was
open until the work of Hess-Schick \cite{hessschick98} who found a
trick in order to control the sign of $\log{\det}_\Gamma(\Delta_p^c)$ in terms of $p$ and $n$. So they were able to show that there is a {\em positive} rational number $q_n$ such that
\begin{equation}\label{lovely} 
T^{(2)}(H^n)=\left(-\frac{1}{\pi}\right)^\frac{n-1}{2} q_n \ .
\end{equation}
Here the metric on $H^n$ is normalized to have sectional curvature $-1$. 
Along the same lines Hess \cite{hessdiss} obtained
\begin{equation}\label{smith}
 T^{(2)}(SL(3,\R)/SO(3))\ne 0\ .
\end{equation}
Let us introduce
$$ Q_n=\frac{2q_n}{(\frac{n-1}{2})!}$$ 
and rewrite (\ref{lovely}) as
\begin{equation}\label{love} 
T^{(2)}(H^n)=(-1)^\frac{n-1}{2}\frac{\pi Q_n}{\vol(S^n)} \ .
\end{equation}
Recall that $S^n$ is the compact dual of $H^n$. The rational number $Q_n$ which does not depend on the normalization of the metric has a nice 
interpretation in terms of Weyl's dimension polynomial for finite-dimensional
representations, see Proposition \ref{mary}. But its significance remains to be clarified further, and it seems to be difficult to write
down a practical formula valid for all odd $n$. One has $Q_3=\frac{1}{3}$, 
$Q_5=\frac{31}{45}$, $Q_7=\frac{221}{210}$. For further information see
\cite{hessschick98}.  

We will reduce Theorem \ref{main} (c) to (\ref{smith}) and the positivity of $Q_n$. Let $X$ be an arbitrary symmetric space satisfying $m(X)=1$. By
the classification of simple Lie groups
$X=X_1\times X_0$, where $m(X_0)=0$ and $X_1=SL(3,\R)/SO(3)$ or $X_1=X_{p,q}:=
SO(p,q)^0/SO(p)\times SO(q)$ for $p,q$ odd. (Here as throughout the paper an upper subscript $0$ denotes the connected component of the identity.) Note that a corresponding decomposition of $Y$ does not necessarily exist. We show

\begin{prop}\label{duality}\mbox{ }\\
\begin{enumerate}
\item[(a)] $T^{(2)}(X_{p,q})=\D(-1)^\frac{pq-1}{2}\chi(X_{p-1,q-1}^d)\frac{\pi Q_{p+q-1}}{\vol(X_{p,q}^d)}.$\\ 
\item[(b)] If $m(X)=1$, then $ T^{(2)}(X)=\D\frac{(-1)^\frac{n_0}{2}\chi(X^d_0)}{\vol(X^d_0)}T^{(2)}(X_1)$.
Here $n_0=\dim X_0$.
\end{enumerate}
\end{prop}
Note that $X_{p-1,q-1}^d=SO(p+q-2)/SO(p-1)\times SO(q-1)$ and for $p,q>1$ $$\chi(X_{p-1,q-1}^d)=2{\frac{p+q-2}{2}\choose \frac{p-1}{2}}\ .$$ 
In fact, it is a classical result that $\chi(X^d)>0$ whenever $m(X)=0$ (compare Theorem \ref{main} (a)). It is equal to the quotient of the orders of
certain Weyl groups (see Section \ref{to}). Now Theorem \ref{main} (c) follows from Proposition \ref{duality}, (\ref{smith}) and the positivity of $Q_n$. 

We are also able to identify the missing constant in (\ref{smith}). 

\begin{prop}\label{hesse}
If $X=SL(3,\R)/SO(3)$, then $\D T^{(2)}(X)=\frac{2}{3}\frac{\pi}{\vol(X^d)}$. 
If the invariant metric on $X$ is induced from twice the trace form
of the standard representation of ${\frak sl}(3,\R)$, then $\vol(X^d)=4\pi^3$, and we have
$$T^{(2)}(X)=\frac{1}{6\pi^2}\ .$$
\end{prop}
In particular, we see that $(-1)^\frac{n-1}{2}T^{(2)}(X)$ is positive for
all $X$ with $m(X)=1$. Proposition \ref{mary} provides a uniform formula
for the $L^2$-torsion of all these spaces. 

\noindent
{\it Acknowledgements}: I am grateful to Wolfgang L\"uck and Thomas Schick
for inspiring discussions which have provided me with a sufficient amount of
motivation and of knowledge on $L^2$-invariants in order to perform the computations which led to the results of the present paper. I am also indebted to Wolfgang L\"uck
for giving me the opportunity to report on them at the Oberwolfach conference ``$L^2$-methods and $K$-theory'', 
September 1999.
In addition, I benefited from discussions with J. Lott, P. Pansu, E. Hess and
U. Bunke.

\section{The Harish-Chandra Plancherel Theorem}\label{shc}

We want to understand the action of the Laplacian and of the corresponding
heat kernels on $L^2(X,\Lambda^pT^*X)$. Since the Laplacian coincides (up to the sign) with
the action of the Casimir operator $\Omega$ of $G$ (Kuga's Lemma \cite{borelwallach80},
Thm. 2.5.) it is certainly enough to understand the "decomposition" of $L^2(X,\Lambda^pT^*X)$ into irreducible unitary representations of $G$. There is an isomorphism
of homogeneous vector bundles $\Lambda^pT^*X\cong G\times_K\Lambda^p\paaa^*$, 
and, hence, of $G$-representations
$$L^2(X,\Lambda^pT^*X)\cong [L^2(G)\otimes \Lambda^p\paaa^*]^K\ .$$
Thus our task consists of two steps: First to understand $L^2(G)$ as a
representation of $G\times G$ which is
accomplished by the Harish-Chandra Plancherel Theorem which we recall in the
present section and, second, to understand spaces of the form $[V_\pi\otimes \Lambda^p\paaa^*]^K$, where $(\pi,V_\pi)$ is an irreducible unitary representation of $G$ which occurs in the Plancherel decomposition. For general $G$,
the second step will resist a naive approach. However, if $\pi(\Omega)=0$, then
the space $[V_\pi\otimes \Lambda^p\paaa^*]^K$ has cohomological meaning, and
the theory of relative $(\gaaa,K)$-cohomology as recalled in the next section
will provide a sufficient amount of information.

An irreducible unitary representation $(\pi, V_\pi)$ of $G$ is called a representation 
of the discrete series if there is a $G$-invariant embedding $V_\pi\hookrightarrow L^2(G)$. Let $\hat G_d$ denote the set of equivalence classes of discrete series representations of $G$. Then we have

\begin{theorem}[Harish-Chandra \cite{HC66}]\label{comp}
$\hat G_d$ is non-empty if and only if $m(X)=0$.
\end{theorem}

Note that $m(X)=0$ means that $G$ has a compact Cartan subgroup. 
The Plancherel Theorem provides a decomposition of $L^2(G)$ which is indexed
by discrete series representations of certain subgroups $M\subset G$ which
we are going to define now. 

Let $\aaaa_0\subset\paaa$ be a maximal abelian subspace. It induces a root
space decomposition
$$ \gaaa=\gaaa_0\oplus \bigoplus_{\alpha\in \Delta(\gaaa,\aaaa_0)} \gaaa_\alpha\ .$$
Choose a decompostion $\Delta(\gaaa,\aaaa_0)=\Delta^+\stackrel{\cdot}{\cup}
-\Delta^+$ into positive and negative roots, and let $\Pi\subset \Delta^+$ be
the subset of simple roots. For any subset $F\subset \Pi$ we define
\begin{eqnarray*} 
\aaaa_F=\{H\in\aaaa_0\:|\: \alpha(H)=0\mbox{ for all }\alpha\in\Pi\}\ ,&\quad&
A_F=\exp(\aaaa_F)\ ,\\
\naaa_F=\bigoplus_{\{\alpha\in \Delta^+| \alpha_{|\aaaa_F}\ne 0\}} 
\gaaa_\alpha\ ,
&\quad& N_F=\exp(\naaa_F)\ .
\end{eqnarray*}
Furthermore, there is a unique (possibly disconnected) subgroup $M_F\subset G$ with Lie algebra $\maaa_F$ such that $M_FA_F$ is the centralizer of $a_F$ in
$G$ and $\maaa_F$ is orthogonal to $\aaaa_F$ with respect to any invariant bilinear form on $\gaaa$. $M_F$ is a reductive subgroup with compact center.
The corresponding parabolic subgroup $P_F:=M_FA_FN_F$ is called a standard parabolic. $P_F$ is called cuspidal if $M_F$ has a compact Cartan subgroup. 
If for two subsets $F,I\subset \Pi$ the spaces $\aaaa_F$ and $\aaaa_I$ are conjugated by an element of $K$ (thus by an element of the Weyl group $W(\gaaa,\aaaa_0)$) we call $P_F$ and $P_I$ associate. (In many cases this already implies $F=I$.) The assignment $P_F\mapsto A_FT$, where $T$ is a compact Cartan subgroup of $M_F$, gives a one-to-one correspondence between association classes of cuspidal parabolic subgroups and conjugacy classes of Cartan subgroups of $G$.

For two subsets $F\subset I\subset\Pi$ we have $P_F\subset P_I$, $M_F\subset M_I$, $A_F\supset A_I$, $N_F\supset N_I$.  If $F$ is understood we will often suppress
the subscript $F$.

For illustration let us consider the two extreme cases. The minimal parabolic arises for $F=\emptyset$. Since $M_\emptyset$ is compact $P_\emptyset$ is always cuspidal.
For $F=\Pi$ we have $P=M=G$, and $G$ is cuspidal iff $m(X)=0$. For any cuspidal
parabolic subgroup we have $\dim A\ge m(X)$, and
there is exactly one association class of cuspidal parabolic subgroups, called
fundamental, with $\dim A=m(X)$. 

Theorem \ref{comp} also holds in the context of such reductive groups like $M$. Thus a parabolic $P=MAN$ is cuspidal iff $M$ has a non-empty discrete series $\hat M_d$.
Let $\aca$ be the complexified dual of the Lie algebra $\aaaa$ of $A$. For a discrete series representation $(\xi,W_\xi)$ of $M$ and $\nu\in\aca$ we form the induced representation
$(\pi_{\xi,\nu}, H^{\xi,\nu})$ by
$$ H^{\xi,\nu}=\left\{ f:G\rightarrow W_\xi\:|\:\begin{array}{c} f(gman)=a^{-(\nu+\rho_\aaaa)}\xi(m)^{-1}f(g)\mbox{ for all }\\ g\in G, man\in MAN,\  f_{|K}\in L^2(K,W_\xi)\end{array}\right\},\ 
(\pi_{\xi,\nu}(g)f)(x)=f(g^{-1}x)\ .$$
Here $\rho_\aaaa=\D\frac{1}{2}\sum_{\alpha\in\Delta^+}\alpha_{|\aaaa}$. If $\nu\in
i\aaaa^*$, then $\pi_{\xi,\nu}$ is unitary. An invariant bilinear form on $\gaaa$ induces corresponding forms on $\maaa$ and $\aca$ and determines Casimir operators $\Omega$ and $\Omega_M$ of $G$ and $M$, respectively. Then we have
\begin{equation}\label{magic}
\pi_{\xi,\nu}(\Omega)=\langle \nu,\nu \rangle -\langle \rho_\aaaa,\rho_\aaaa \rangle+
\xi(\Omega_M)\ .
\end{equation}     
Note that $\xi(\Omega_M)$ is a non-negative real scalar.

Let $\cC(G)\subset L^2(G)$ be the Harish-Chandra Schwartz space (for a definition see e.g. \cite{wallach88}, 7.1.2). It is stable under the left and right regular representions $l$ and $r$ of $G$. Let $\cC(G)_{K\times K}=\{f\in\cC(G)\:|\:
\dim\mbox{span}\{l_{k_1}r_{k_2}f\:|\:k_1,k_2\in K\}<\infty\}$ be the subspace of Schwartz functions which are $K$-finite from the left and the right. Note that $\cC(G)_{K\times K}$ is 
dense in $L^2(G)$. Suppose that $\nu\in
i\aaaa^*$. Then for $f\in\cC(G)$ 
$$ \pi_{\xi,\nu}(f):=\int_G f(g)\pi_{\xi,\nu}(g)\: dg$$
is a well-defined trace class operator on $H^{\xi,\nu}$ which has finite rank if $f\in\cC(G)_{K\times K}$. Note that the map $f\mapsto \pi_{\xi,\nu}(f)$ intertwines
the $G$-actions in the following way: $\pi_{\xi,\nu}(l_x r_yf)=\pi_{\xi,\nu}(x)\pi_{\xi,\nu}(f)\pi_{\xi,\nu}(y^{-1})$ for $x,y\in G$. 

The Harish-Chandra Plancherel Theorem can now be formulated as follows:

\begin{theorem}[Harish-Chandra \cite{HC76III}]\label{HCPl}
For each cuspidal parabolic subgroup as constructed above and any discrete
series representation $\xi$ of the corresponding group $M$ there exists an
explicitly computable analytic function 
$ p_\xi: i\aaaa^*\rightarrow [0,\infty)$ of polynomial growth (the Plancherel density) such that
for any $f\in\cC(G)_{K\times K}$ and $g\in G$
$$ f(g)=\sum_P\sum_{\xi\in\hat M_d} \int_{\aaaa^*} \Tr(\pi_{\xi,i\nu}(f)\pi_{\xi,i\nu}(g^{-1}))\: p_\xi(i\nu)d\nu\ .$$
Here the first sum runs over a set of representatives $P=P_F$ of association classes of cuspidal parabolic subgroups of $G$.
\end{theorem}

For more details on the Plancherel Theorem and the structure theory behind it the interested reader may consult the textbooks \cite{knapp86}, \cite{wallach88}, \cite{wallach92}.

Note that the Plancherel measures $p_\xi(i\nu)d\nu$ depend on the normalization of
the Haar measure $dg$. In the remainder of the paper we use the following one.
Let $dx$ be the Riemannian volume form of $X=G/K$ and $dk$ be the Haar measure of $K$
with total mass one. Then $\D \int_G f(g)\:dg=\int_X \tilde f(x)\:dx$, where $\D \tilde f(gK)=\int_K f(gk)\:dk$. We normalize the invariant bilinear form on
$\gaaa$ such that its restriction to $\paaa\cong T_{eK} X$ coincides with the
Riemannian metric of $X$. Let $d\nu$ be the Lebesgue measure corresponding to the induced form on $\aaaa^*$. By these choices $p_\xi$ is uniquely determined.

We are now able to give a kind of spectral expansion of $\Tr_\Gamma e^{-t\Delta_p^*}$.

\begin{kor}\label{huhu}
\begin{eqnarray} 
\Tr_\Gamma e^{-t\Delta_p}&=&\vol(Y)\sum_P\sum_{\xi\in\hat M_d} \int_{\aaaa^*} e^{-t(\|\nu\|^2+\|\rho_\aaaa\|^2-\xi(\Omega_M))}\dim [ H^{\xi,i\nu}\otimes \Lambda^p\paaa^*]^K\: p_\xi(i\nu)d\nu\label{ms1}\\
&=& \vol(Y)\sum_P\sum_{\xi\in\hat M_d} \int_{\aaaa^*} e^{-t(\|\nu\|^2+\|\rho_\aaaa\|^2-\xi(\Omega_M))}\dim [ W_\xi\otimes \Lambda^p\paaa^*]^{K_M}\: p_\xi(i\nu)d\nu\ .\label{ms2}
\end{eqnarray}
Here $K_M:=K\cap M$ denotes the maximal compact subgroup of $M$. There are only finitely many pairs $(P,\xi)$ with $[ H^{\xi,i\nu}\otimes \Lambda^p\paaa^*]^K\cong [ W_\xi\otimes \Lambda^p\paaa^*]^{K_M}\ne\{0\}$.
\end{kor}
\proof We define $k_t\in[\cC(G)\otimes\End(\Lambda^p\paaa^*)]^{K\times K}$
by $k_t(g):= e^{-t\Delta_p}(eK,gK)\circ g$. For $f\in
L^2(X,\Lambda^pT^*X)\cong [L^2(G)\otimes \Lambda^p\paaa^*]^K$ we have
$e^{-t\Delta_p}f(g_0)=\D\int_G k_t(g) f(g_0g)\:dg$. In addition, $\tr\, e^{-t\Delta_p}(x,x)=\tr\, k_t(e)$ for any $x\in X$. We consider $\pi_{\xi,i\nu}(k_t)$ as an operator acting on $H^{\xi,i\nu}\otimes \Lambda^p\paaa^*$. Using Kuga's Lemma, Equation (\ref{magic}), and the $K\times K$-invariance of $k_t$ one derives that $\pi_{\xi,i\nu}(k_t)=e^{-t(\|v\|^2+\|\rho_\aaaa\|^2-\xi(\Omega_M))}P$, where $P$ is
the orthogonal projection onto the subspace of $K$-invariants in $H^{\xi,i\nu}\otimes \Lambda^p\paaa^*$. The Plancherel formula now yields
$$\tr\, k_t(e)=\sum_P\sum_{\xi\in\hat M_d} \int_{\aaaa^*} e^{-t(\|v\|^2+\|\rho_\aaaa\|^2-\xi(\Omega_M))}\dim [ H^{\xi,i\nu}\otimes \Lambda^p\paaa^*]^K\: p_\xi(i\nu)d\nu\ .$$
This proves (\ref{ms1}).
Since $H^{\xi,i\nu}\cong L^2(K\times_{K_M} W_\xi)$ as a representation of $K$
Equation (\ref{ms2}) follows by Frobenius reciprocity. The last assertion is
a consequence of the Blattner formula (see e.g. \cite{wallach88}, 6.5.4) for
the $K_M$-types of discrete series representations of $M$.
\hB

\section{$({\frak g}, K)$-cohomology}\label{coho}

If $(\pi, V_\pi)$ is a representation of $G$ on a complete locally convex
Hausdorff topological vector space, then we can form its subspace $V_{\pi,K}$
consisting of all $K$-finite smooth vectors of $V_\pi$. $V_{\pi,K}$ becomes a simultaneous module under $\gaaa$ and $K$, where both actions satisfy
the obvious compatibility conditions. Such a module is called a $(\gaaa,K)$-module (see \cite{borelwallach80}, 0.2). 

We are interested
in the functor of $(\gaaa,K)$-cohomology $V\mapsto H^*(\gaaa,K,V)$ which
goes from the category of $(\gaaa,K)$-modules to the category of vector spaces.
It is the right derived functor of the left exact functor taking $(\gaaa,K)$-invariants. $H^*(\gaaa,K,V)$ can be computed using the standard relative
Lie algebra cohomology complex $([V\otimes \Lambda^*\paaa^*]^K, d)$, where
$$ d\omega(X_0,\dots,X_p)=\sum_{i=0}^{p}(-1)^i \pi(X_i)\omega(X_0,\dots,\hat X_i,\dots, X_p)\ ,\quad \omega\in [V\otimes \Lambda^p\paaa^*]^K,\  X_i\in\paaa\ .$$
Note that for $V=C^\infty(G)_K$ this complex is isomorphic to the de Rham
complex of the symmetric space $X$.

Let $\cZ(\gaaa)$ be the center of the universal enveloping algebra of $\gaaa$. 
If $(\pi,V)$ is an irreducible $(\gaaa,K)$-module, then any $z\in \cZ(\gaaa)$
acts by a scalar $\chi_{\pi}(z)$ on $V$. The homomorphism 
$\chi_\pi:\cZ(\gaaa)\rightarrow \C$ is called the infinitesimal character of $V$. The following basic result can be considered as an algebraic version of Hodge theory. 

\begin{prop}[\cite{borelwallach80}, II.3.1. and I.5.3.]\label{hodge}
Let $(\pi,V_\pi)$ be an irreducible unitary representation of $G$, and let 
$(\tau,F)$ be an irreducible finite-dimensional representation of $G$. Then
$$ H^p(\gaaa,K,V_{\pi,K}\otimes F)=\left\{
\begin{array}{cc} 
[V_{\pi}\otimes F\otimes \Lambda^p\paaa^*]^K&\pi(\Omega)=\tau(\Omega)\\
\{0\}&\pi(\Omega)\ne \tau(\Omega)
\end{array}\right. \ .$$
If $H^p(\gaaa,K,V_{\pi,K}\otimes F)\ne\{0\}$, then $\chi_\pi=\chi_{\tilde \tau}$, where $\tilde\tau$ is the dual representation of $\tau$.
\end{prop}

The cohomology groups $H^p(\gaaa,K,V_{\pi,K}\otimes F)$ for the representations $\pi=\pi_{\xi,i\nu}$ occuring in the Plancherel Theorem have been computed in \cite{borelwallach80}. We shall need
the following information.

\begin{prop}[\cite{borelwallach80},II.5.3., III.5.1., and III.3.3.]\label{bobo}\mbox{}\\
\begin{enumerate}
\item[(a)] Let 
$(\tau,F)$ be an irreducible finite-dimensional representation of $G$ and $\pi\in\hat G_d$ with $\chi_\pi=\chi_{\tilde \tau}$. Then
$$ \dim H^p(\gaaa,K,V_{\pi,K}\otimes F)=\left\{
\begin{array}{cc} 
1&p=\frac{n}{2}\\
0&\mbox{otherwise}
\end{array}\right. \ .$$
\item[(b)] Let $(\pi_{\xi,\nu}, H^{\xi,\nu})$ a representation occuring in the
Plancherel Theorem. Then 
$$H^*(\gaaa,K,H^{\xi,i\nu}_K)=\{0\}$$
unless $P$ is fundamental, $\nu=0$ and $\xi$ belongs to a certain non-empty finite
subset $\Xi\subset\hat M_d$. If $P$ is fundamental and $\xi\in\Xi$, then
$$ \dim H^p(\gaaa,K,H^{\xi,0}_K)=\left\{
\begin{array}{cc} 
{m\choose p-\frac{n-m}{2}}&p\in [\frac{n-m}{2},\frac{n+m}{2}]\\
0&\mbox{otherwise}
\end{array}\right. \ .$$
\end{enumerate}
\end{prop}

Choose a Cartan subalgebra $\haaa\subset \gaaa$ and a system of positive roots. Via the Harish-Chandra isomorphism any infinitesimal character $\chi_\pi: \cZ(\gaaa)\rightarrow \C$ is given by an element $\Lambda_\pi\in\haaa^*_\C$, which is uniquely determined up to the action of the Weyl group $W(\gaaa,\haaa)$ (\cite{wallach88}, 3.2.4.). The set of infinitesimal characters of discrete series
representations coincides (in case $m=0$) with the set of infinitesimal characters of finite-dimensional representations, which are of the form $\mu_\tau+\rho_\gaaa$, $\mu_\tau$ and $\rho_\gaaa$ being the highest weight of $\tau$ and the half-sum of positive roots, respectively. In the following we will represent infinitesimal characters of discrete series representations by elements of this form. We introduce a partial
order on $\haaa^*_\C$ by saying that $\mu>\nu$, if $\mu-\nu$ is a sum of (not necessarily distinct) positive roots. Then a careful examination of the proof of Proposition II.5.3.
in \cite{borelwallach80}, which only rests on some basic knowledge of the possible $K$-types occuring in discrete series representations, shows that slightly more
than Proposition \ref{bobo} (a) is true.

\begin{prop}\label{obob}
Let 
$(\tau,F)$ be an irreducible finite-dimensional representation of $G$ and $\pi\in\hat G_d$ with $\Lambda_\pi\not<\Lambda_{\tilde \tau}$. Then
$$ \dim [V_{\pi,K}\otimes F\otimes \Lambda^p\paaa^*]^K=\left\{
\begin{array}{cc} 
1&p=\frac{n}{2}, \chi_\pi=\chi_{\tilde \tau}\\
0&\mbox{otherwise}
\end{array}\right. \ .$$
\end{prop}

\section{$L^2$-Betti numbers and Novikov-Shubin invariants}\label{co}

In this section we shall prove parts (a) and (b) of Theorem \ref{main} as well
as Proposition \ref{spec}.

Let $L^2(X,\Lambda^pT^*X)_d$ be the discrete subspace of $L^2(X,\Lambda^pT^*X)$,
i.e., the direct sum of the $L^2$-eigenspaces of the Laplacian. The Plancherel Theorem in particular says that as a representation of $G$
$$L^2(X,\Lambda^pT^*X)_d \cong \bigoplus_{\pi\in \hat G_d} 
V_{\tilde \pi}\otimes 
[V_\pi\otimes\Lambda^p\paaa^*]^K\ .$$
Let $(\tau_0, \C)$ be the trivial representation of $G$. Let $\pi\in \hat G_d$. Then $\Lambda_\pi\not<\Lambda_{\tau_0}=\rho_\gaaa$. Proposition \ref{obob} now yields
\begin{equation}\label{flu} \dim[V_\pi\otimes\Lambda^p\paaa^*]^K=\left\{\begin{array}{cc}
1& p=\frac{n}{2}, \chi_\pi=\chi_{\tau_0}\\ 0&\mbox{otherwise}\end{array}\right.\ .
\end{equation}
Since in case $m=0$ discrete series representations with infinitesimal character $\chi_{\tau_0}$ always
exist (see e.g. \cite{knapp86}, Thm. 9.20 or \cite{wallach88}, Thm. 6.8.2) and  $\tau_0(\Omega)=0$ this implies Proposition \ref{spec}.
In particular, $b_p^{(2)}(Y)$
is non-zero exactly when $m=0$ and $p=\frac{n}{2}$. Using that the $L^2$-Euler characteristic
coincides with the usual Euler characteristic and applying Hirzebruch proportionality we obtain 
\begin{equation}\label{bor}
b_p^{(2)}(Y)=(-1)^\frac{n}{2}\chi(Y)=\D\frac{\vol(Y)}{\vol(X^d)}\chi(X^d)\ .
\end{equation} 
We will give an alternative, purely analytic proof of that formula in the sequel
of Corollary \ref{jandl}. This finishes the proof of part (a) of Theorem \ref{main}.

We now turn to part (b). In order to compute $\alpha_p(Y)$ we need an
expression for $\Tr_\Gamma e^{-t\Delta_p^c}$.

\begin{prop}\label{esmuss}
For any triple $(P,\xi,\nu)$ appearing in (\ref{ms1}) let 
$B^p(\xi,\nu)=d([ H^{\xi,i\nu}\otimes \Lambda^{p-1}\paaa^*]^K)$ be the space
of coboundaries in the relative Lie algebra cohomology complex and $b^p(\xi,\nu)$ be its dimension. Then 
\begin{equation}\label{ms3}
\Tr_\Gamma e^{-t\Delta^c_p}=\vol(Y)\sum_{P\ne G}\sum_{\xi\in\hat M_d} \int_{\aaaa^*} e^{-t(\|\nu\|^2+\|\rho_\aaaa\|^2-\xi(\Omega_M))}b^{p+1}(\xi,\nu)\: p_\xi(i\nu)d\nu\ .
\end{equation}
Here the first sum runs over a set of representatives $P=P_F$ of association classes of proper cuspidal parabolic subgroups of $G$. If $m>0$, $P$ is fundamental, $\xi\in\Xi$ (see Proposition \ref{bobo} (b)), and $\nu\ne 0$, then 
\begin{equation}\label{ms4}
b^{p+1}(\xi,\nu)=\left\{
\begin{array}{cc} 
{m-1\choose p-\frac{n-m}{2}}&p\in [\frac{n-m}{2},\frac{n+m}{2}-1]\\
0&\mbox{otherwise}
\end{array}\right. \ .
\end{equation}
\end{prop}
\proof We proceed exactly as in the proof of Corollary \ref{huhu}. The
kernel $e^{-t\Delta^c_p}$ determines a function $k_t^c\in[\cC(G)\otimes\End(\Lambda^p\paaa^*)]^{K\times K}$. Then one
computes that $\pi_{\xi,i\nu}(k_t^c)=e^{-t(\|v\|^2+\|\rho_\aaaa\|^2-\xi(\Omega_M))}P^c$,
where $P^c$ is the projection to the orthogonal complement in 
$[H^{\xi,i\nu}\otimes \Lambda^{p}\paaa^*]^K$ of the space of cocycles
in the relative Lie algebra cohomology complex. The dimension of that 
complement is equal to $b^{p+1}(\xi,\nu)$. If $m=0$, $P=G$, and $\xi\in \hat G_d$, then $b^{p+1}(\xi)=0$ for all $p$ by (\ref{flu}). This proves (\ref{ms3}).

Let now $m>0$, $P$ be fundamental, $\xi\in\Xi$, and $\nu\in\aaaa^*$. In this case $H^p(\gaaa,K,H^{\xi,0}_K)\ne\{0\}$ for some $p$. Hence Proposition \ref{hodge} implies that $\pi_{\xi,0}(\Omega)=0$ and  
$$h^p(\xi):=\dim [H^{\xi,i\nu}\otimes \Lambda^{p}\paaa^*]^K=\dim [H^{\xi,0}\otimes \Lambda^{p}\paaa^*]^K=\dim H^p(\gaaa,K,H^{\xi,0}_K)\ .$$
By Proposition \ref{bobo} (b) we have 
$ h^p(\xi)={m\choose p-\frac{n-m}{2}}$.
On the other hand, $\dim H^p(\gaaa,K,H^{\xi,i\nu}_K)=0$
for $\nu\ne 0$ implies that $h^p(\xi)=b^p(\xi,\nu)+b^{p+1}(\xi,\nu)$. (\ref{ms4}) now follows inductively.
\hB

If $H^p(\gaaa,K,H^{\xi,i\nu}_K)=\{0\}$ for all
$\nu\in\aaaa^*$, then by Proposition \ref{hodge} $\dim [ H^{\xi,i\nu}\otimes \Lambda^p\paaa^*]^K=0$ (which is independent of $\nu$) or $\D \inf_{\nu\in\aaaa^*}(\|\nu\|^2+\|\rho_\aaaa\|^2-\xi(\Omega_M))>0$.
Now Proposition \ref{esmuss} implies that the spectrum of $\Delta^c_p$ has
a gap around zero, which means $\alpha_p(Y)=\infty^+$, unless $m>0$ and $p\in [\frac{n-m}{2},\frac{n+m}{2}-1]$.
In the latter case we obtain for some $c>0$
$$\Tr_\Gamma e^{-t\Delta^c_p}=\vol(Y){m-1\choose p-\frac{n-m}{2}}\sum_{\xi\in\Xi}
\int_{\aaaa^*} e^{-t\|\nu\|^2}\: p_\xi(i\nu)d\nu+ O(e^{-ct})\ \mbox{ as } t\to\infty\ ,$$
where $\aaaa$ corresponds to a fundamental parabolic subgroup and thus has
dimension $m$.

Let $P=MAN$ be fundamental and $\xi\in\hat M_d$. Then $\dim \naaa=:2u$ is even.  Choose a compact Cartan subgroup $T\subset M$ with Lie algebra $\taaa$ and a system $\Delta^+(\maaa,\taaa)$ of positive roots.
Considering the pair $(\maaa,\taaa)$ instead of $(\gaaa,\haaa)$ we can define
$\Lambda_\xi\in\taaa_\C^*$ in the same way as at the end of Section \ref{coho}. 
Then $\haaa:=\taaa\oplus \aaaa$ is a Cartan subalgebra of $\gaaa$. Let $\Phi^+$
be a system of positive roots for $(\gaaa,\haaa)$ containing $\Delta^+(\maaa,\taaa)$. Then there exists a positive constant $c_X$ depending
only on the normalization of the volume form $dx$ such that
\begin{equation}\label{pudels}
p_\xi(\nu)=c_X (-1)^u\prod_{\alpha\in\Phi^+}
\frac{\langle \alpha,\Lambda_\xi+\nu\rangle}{\langle \alpha,\rho_\gaaa\rangle}
\end{equation}
(see \cite{HC76III}, Thm. 24.1, \cite{wallach92}, Thm. 13.5.1 or \cite{knapp86},
Thm. 13.11). In particular, $p_\xi$ is an even polynomial of degree $\dim\naaa$.
The factor $(-1)^u$ makes it nonnegative on $i\aaaa^*$.

The element $\Lambda_\xi+\nu$ gives the infinitesimal character of $\pi_{\xi,\nu}$: $\Lambda_{\pi_{\xi,\nu}}=\Lambda_\xi+\nu$ (\cite{knapp86}, Prop. 8.22).
Let now $\xi\in\Xi$. Then Propositions \ref{bobo} (b) and \ref{hodge} imply
that $\pi_{\xi,0}$ has the same infinitesimal character as the trivial representation. It follows that $\Lambda_\xi$ is conjugated in $\haaa_\C^*$ by an element of the Weyl
group $W(\gaaa,\haaa)$ to $\rho_\gaaa$. By (\ref{pudels}) we obtain 
$p_\xi(0)=\pm c_X\ne 0$. On the other hand $p_\xi(0)\ge 0$, hence 
\begin{equation}\label{crux}
p_\xi(0)>0\ .
\end{equation} 
We decompose $p_\xi(\nu)=\D\sum_{k=0}^{u} p_{\xi,2k}(\nu)$ into homogeneous polynomials.
Set $q_{\xi,k}=\D\int_{\|\nu\|=1} p_{\xi,2k}(\nu)\:d\nu$. Then $q_{\xi,0}>0$ and
$$\int_{\aaaa^*} e^{-t\|\nu\|^2}\: p_\xi(i\nu)d\nu=
\sum q_{\xi,k}\int_0^\infty e^{-tr^2} r^{m-1+2k}dr=\sum t^{-(\frac{m}{2}+k)}q_{\xi,k}\int_0^\infty e^{-y^2}y^{m-1+2k}dy\ .$$
Thus for $p\in [\frac{n-m}{2},\frac{n+m}{2}-1]$ the leading term of $\Tr_\Gamma e^{-t\Delta^c_p}$ as $t\to\infty$ is a non-zero multiple of $t^{-\frac{m}{2}}$.
This completes the proof of Theorem \ref{main} (b).

\section{$L^2$-torsion}\label{to}

For even dimensional manifolds the $L^2$-torsion vanishes. Thus we may assume that $m$ is odd, in particular $m\ge 1$. Then $\Delta_p^\prime=\Delta_p$.
We first want to compute
$$ k_X(t):=\frac{1}{2\vol(Y)}\sum_{p=0}^{n} (-1)^pp\,\Tr_\Gamma e^{-t\Delta_p}\ .$$
Then
$$ T^{(2)}(X)=\frac{d}{ds}_{|s=0}\left(\frac{1}{\Gamma(s)}\int_0^\ve
k_X(t)t^{s-1}dt\right)+\int_\ve^\infty k_X(t)t^{-1}dt\ . $$
Let $P=MAN$ be a parabolic subgroup appearing in (\ref{ms2}). Set $K_M=K\cap M$ and $\paaa_\maaa=\paaa\cap \maaa$. Then an elementary calculation in the representation ring
$R(K_M)$ of $K_M$ yields
$$
\sum_{p=0}^{n} (-1)^pp\,\Lambda^p\paaa^*=0\ ,\qquad \mbox {if }\dim\aaaa\ge 2\ ,$$
and
\begin{eqnarray*}
\sum_{p=0}^{n} (-1)^pp\,\Lambda^p\paaa^*&=&\sum_{p=0}^{n-1} (-1)^{p+1}\Lambda^p(\paaa_\maaa^*\oplus \naaa^*)\\
&=&\sum_{l=0}^{\dim\naaa} (-1)^{l+1}(\Lambda^{ev}\paaa_\maaa^*-\Lambda^{odd}\paaa_\maaa^*)
\otimes\Lambda^l\naaa\ ,\quad\mbox{if }\dim\aaaa=1
\end{eqnarray*}
(see \cite{moscovicistanton91}, Prop. 2.1 and Lemma 2.3). It follows
from (\ref{ms2}) that $k_X(t)\equiv 0$ for $m>1$, hence $\rho^{(2)}(Y)=0$. 

From now on let $m=1$.
Let $P=MAN$ be a fundamental parabolic subgroup of $G$. Then (\ref{ms2}) gives
\begin{eqnarray}
k_X(t)&=&\frac{1}{2} \sum_{l=0}^{\dim\naaa} (-1)^{l+1}\sum_{\xi\in\hat M_d}
\dim\, [W_\xi\otimes(\Lambda^{ev}\paaa^*_\maaa-\Lambda^{odd}\paaa^*_\maaa)
\otimes\Lambda^l\naaa^*]^{K_M}
\nonumber\\&&
\hspace{5cm}\int_{\aaaa^*}e^{-t(\|\nu\|^2+\|\rho_\aaaa\|^2-\xi(\Omega_M))}
\: p_\xi(i\nu)d\nu\ . \label{aha}
\end{eqnarray}

Now $X=X_1\times X_0$, $X_1=G_1/K_1$, $X_0=G_0/K_0$, $m(X_0)=0$ as explained in the introduction. Although (\ref{aha}) can be evaluated directly for general
$X$ with $m(X)=1$ we prefer to reduce the computation to the irreducible case
$X=X_1$. In order to compute $T^{(2)}(X)$ it is sufficient to compare $\rho^{(2)}(Y)$ with $\vol(Y)$ for one particular $Y=\Gamma\backslash X$.
If we choose $\Gamma$ of the form $\Gamma_1\times \Gamma_0$, where  $\Gamma_0\subset G_0$ and
$\Gamma_1\subset G_1$, then 
$$\rho^{(2)}(Y)= \chi(Y_0)\rho^{(2)}(Y_1)\ , $$
where $Y_1=\Gamma_1\backslash X_1$, $Y_0=\Gamma_0\backslash X_0$.
Applying Hirzebruch proportionality we obtain the assertion of Proposition \ref{duality} (b)
\begin{eqnarray} 
T^{(2)}(X)&=&\frac{\rho^{(2)}(Y)}{\vol(Y_0)\vol(Y_1)}
=\frac{\chi(Y_0)}{\vol(Y_0)}T^{(2)}(X_1)\nonumber\\
&=&\frac{(-1)^\frac{n_0}{2}\chi(X^d_0)}{\vol(X^d_0)}T^{(2)}(X_1)\ .\label{hicks}
\end{eqnarray}

It remains to deal with the case $X=X_1$. We can assume that $G=SO(p,q)^0$, $p\le q$ odd, or $G=SL(3,\R)$. Then $M\cong SO(p-1,q-1)$, $\naaa\cong \R^{p+q-2}$ or $M\cong {}^0GL(2,\R):=\{A\in GL(2,\R)\:|\:|\det A|=1\}$,   
$\naaa\cong \R^{2}$, respectively, and $M$ acts on $\naaa$ via the standard representation. Note that $M$ is not connected unless $G=SO(1,q)^0$.
The $M^0$-representations $\Lambda^l\naaa^*\otimes\C$ are irreducible
unless $G=SO(p,q)$ and $l=u=\frac{1}{2}\dim\naaa$. In the latter case $\Lambda^u\naaa^*\otimes\C$
decomposes into two irreducible components $\Lambda^+\naaa\oplus\Lambda^-\naaa$. Since compact Cartan subgroups of $M$ are connected the discrete series representations of $M$ are induced from discrete
series representations of $M^0$: $W_\xi=\mbox{Ind}_{M^0}^{M}(W_{\xi^0})$,
$\xi_0\in(\hat M^0)_d$ (see \cite{wallach88}, 6.9 and 8.7.1). As representations of $K_M$ we have $W_\xi\cong\mbox{Ind}_{K_M^0}^{K_M}(W_{\xi^0})$. By Frobenius reciprocity we obtain
$$\dim [W_\xi\otimes(\Lambda^{ev}\paaa_\maaa^*-\Lambda^{odd}\paaa_\maaa^*)
\otimes\Lambda^l\naaa^*]^{K_M}=\dim [W_{\xi_0}\otimes(\Lambda^{ev}\paaa_\maaa^*-\Lambda^{odd}\paaa_\maaa^*)
\otimes\Lambda^l\naaa^*]^{K_M^0}\ .$$
Note that the infinitesimal characters $\chi_\xi$ and $\chi_{\xi_0}$ coincide.
By $\chi(\maaa,K^0_M,.)$ we denote the Euler characteristic of relative Lie
algebra cohomology. Set $v=\frac{1}{2}\dim\paaa_\maaa$. Applying Propositions \ref{hodge} and \ref{bobo} (a) to
$M^0$ instead of $G$ we obtain
$$\dim [W_{\xi_0}\otimes(\Lambda^{ev}\paaa_\maaa^*-\Lambda^{odd}\paaa_\maaa^*)
\otimes\Lambda^*\naaa^*]^{K_M^0}=\chi(\maaa,K^0_M,W_{\xi_0,K^0_M}\otimes\Lambda^*\naaa^*)=
\left\{\begin{array}{cc}
(-1)^v& \chi_{\xi}=\chi_{\Lambda^*\naaa}\\ 0&\mbox{otherwise}\end{array}\right. .$$
Here $\Lambda^*\naaa^*$, $*=l,+,-$, denotes an irreducible component of $\Lambda^l\naaa^*\otimes\C$. 

In all cases under consideration the set $\{\alpha_{|\aaaa}\:|\: \alpha\in \Delta^+, \alpha_{|\aaaa}\ne 0 \}$ consists of a single element $\alpha_0\in \aaaa^*$. It follows that $\rho_\aaaa=u\alpha_0$. Moreover, 
$\Omega_M$ acts on $\Lambda^l\naaa$ as $l(2u-l)\|\alpha_0\|^2\id$ (compare \cite{moscovicistanton91}, Lemma 2.5).

In order to evaluate (\ref{aha}) further we have to determine the constant
$c_X$ in formula (\ref{pudels}). This can be done in complete generality.
So for a moment we drop the assumptions $m=1$, $X=X_1$.

\begin{lem}\label{kern}
Let $P=MAN\subset G$ be fundamental. We retain the notation introduced
before (\ref{pudels}). Set $W_A=\{k\in K\:|\:\Ad(k)\aaaa\subset\aaaa\}/K_M$, 
$S^d_A=\exp(i\aaaa)K\subset X^d$, and let $\Phi^+_\kaaa$ be a positive root system for $(\kaaa,\taaa)$ with corresponding half sum $\rho_\kaaa$. Then
\begin{eqnarray}
c_X&=&\frac{1}{|W_A|(2\pi)^\frac{n+m}{2}}
\frac{\prod_{\alpha\in\Phi^+}\langle \alpha,\rho_\gaaa\rangle}{\prod_{\alpha\in\Phi_\kaaa^+}\langle \alpha,\rho_\kaaa\rangle}\label{ottos}\\
&=&\frac{1}{|W_A|}\frac{\vol(S^d_A)}{(2\pi)^m}\frac{1}{\vol(X^d)}\label{mops}\ .
\end{eqnarray}
\end{lem}
\proof Formula (\ref{ottos}) is a combination of \cite{HC75I}, Thm. 37.1, with
\cite{HC76III}, Cor. 23.1, Thm. 24.1 and Thm. 27.3. In order to
apply these results correctly one has to take into account that
Harish-Chandras and our normalizations of the measures $dg$ and $d\nu$ all of them starting from a
fixed invariant bilinear form on $\gaaa$ differ by the factors 
$2^\frac{n-\dim \aaaa_0}{2}$ (\cite{HC75I}, Section 7 and Lemma 37.2) and $(2\pi)^m$, respectively.
On the other hand we have 
\begin{equation}\label{ogottogott}
\prod_{\alpha\in\Phi_\kaaa^+}\langle \alpha,\rho_\kaaa\rangle=(2\pi)^\frac{\dim K/T}{2}\frac{\vol(T)}{\vol(K)}
\end{equation}
(see e.g. \cite{HC75I}, Lemma 37.4). Here the volumes are the Riemannian ones
corresponding to the invariant bilinear form $\langle.,.\rangle$. Formula (\ref{ogottogott}) holds for any pair $(K,T)$ of a connected compact Lie group
and a maximal torus. Applying it also to the pair $(G^d, H^d)$, where $H^d$ is the
maximal torus of $G^d$ with Lie algebra $\taaa\oplus i\aaaa$, we obtain  
\begin{equation}\label{trotz}
\frac{\prod_{\alpha\in\Phi^+}\langle \alpha,\rho_\gaaa\rangle}{\prod_{\alpha\in\Phi_\kaaa^+}\langle \alpha,\rho_\kaaa\rangle}=(2\pi)^\frac{n-m}{2}\frac{\vol(K)\vol(H^d)}
{\vol(G^d)\vol(T)}=(2\pi)^\frac{n-m}{2}\frac{\vol(S^d_A)}{\vol(X^d)}\ .
\end{equation}
The second equality follows from the fact that the map $H^d/T\rightarrow X^d$, 
$hT\mapsto hK$, is an isometric embedding with image $S^d_A$. This proves (\ref{mops}). 
\hB

In particular, specializing (\ref{pudels}) and (\ref{mops}) to the case $m=0$
we obtain as a consequence of Weyl's dimension formula

\begin{kor}\label{jandl}
Let $\pi$ be a discrete series representation of $G$ having the same infinitesimal character as the finite-dimensional representation $\tau$, then
$$p_\pi=\frac{\dim \tau}{\vol(X^d)}\ .$$
\end{kor}
Let us now give the promised analytic proof of (\ref{bor}). For a fixed finite-dimensional representation $\tau$ of $G$ (which still is assumed to be connected) there are exactly 
$|W(\gaaa,\taaa)|/|W(\kaaa,\taaa)|$ equivalence classes of discrete series 
representations with infinitesimal character $\chi_\tau$ (see \cite{wallach88}, Thm. 8.7.1, or \cite{knapp86}, Thm. 12.21). But this quotient of orders of Weyl groups is equal to $\chi(X^d)$ (see e.g. \cite{bott65}). 
By (\ref{ms1}), (\ref{flu}) and Corollary \ref{jandl} we obtain
$$ b_{\frac{n}{2}}^{(2)}=\vol(Y)\sum_{\pi\in\hat G_d, \chi_\pi=\chi_{\tau_0}} p_\pi
=\vol(Y)\chi(X^d)\frac{1}{\vol(X^d)}\ .$$

We return to the evaluation of $k_X(t)$ for $m=1$, $X=X_1$. The polynomial
$p_\xi$ only depends on the infinitesimal character of $\xi$. By $\Lambda_*\in i\taaa^*$, $*=l,+,-$, we denote the infinitesimal character of the irreducible $M$-representation $\Lambda^*\naaa$. We set
$$ p_l(\nu):= \left\{\begin{array}{cc}
\prod_{\alpha\in\Phi^+}
\frac{\langle \alpha,\Lambda_l+\nu\rangle}{\langle \alpha,\rho_\gaaa\rangle}& G=SL(3,\R)\mbox{ or } l\ne u \\ 
\prod_{\alpha\in\Phi^+}\frac{\langle \alpha,\Lambda_++\nu\rangle}{\langle \alpha,\rho_\gaaa\rangle}
+\prod_{\alpha\in\Phi^+}\frac{\langle \alpha,\Lambda_-+\nu\rangle}{\langle \alpha,\rho_\gaaa\rangle}=2\prod_{\alpha\in\Phi^+}\frac{\langle \alpha,\Lambda_++\nu\rangle}{\langle \alpha,\rho_\gaaa\rangle}&\mbox{otherwise}\end{array}\right. .$$
For fixed infinitesimal character there are $|W(\maaa,\taaa)|/|W_{K_M}|$ equivalence classes of discrete series 
representations, where $W_{K_M}=\{k\in K_M\:|\:\Ad(k)\taaa\subset\taaa\}/T$.
Note that there is an embedding $W(\kaaa_\maaa,\taaa)\hookrightarrow W_{K_M}$ which becomes an isomorphism if $K_M$ is connected. Furthermore, $|W_A|=1$ except for $G=SO(1,q)^0$, where $|W_A|=2$. In any
case $|W_{K_M}||W_A|=2|W(\kaaa_\maaa,\taaa)|$. Let $X^d_M$ be the compact
dual of $X_M=M/K_M=M^0/K_M^0$. Then $|W(\maaa,\taaa)|/|W_{K_M}| |W_A|=\frac{1}{2}\chi(X^d_M)$. In addition, $u+v=\frac{n-1}{2}$ and $$\frac{\vol(S^d_A)}{2\pi}=\frac{1}{\|\alpha_0\|}\ .$$
Summarizing the above discussion we obtain 
$$ k_X(t)=(-1)^\frac{n-1}{2}\frac{\chi(X^d_M)}{4\|\alpha_0\|\vol(X^d)}\sum_{l=0}^{2u} (-1)^{l+1} k_l(t)\ ,
$$
where
$$ k_l(t)=\int_{\aaaa^*} e^{-t(\|\nu\|^2+(u-l)^2\|\alpha_0\|^2)}p_l(i\nu)d\nu\ .$$  
For an even polynomial $P$ and $c\ge 0$ set
$$ k_{P,c}(t)=\int_{-\infty}^{\infty} e^{-t(y^2+c^2)}P(iy)dy\ .$$
Then (compare \cite{fried862}, Lemma 2 and Lemma 3, \cite{mathai92}, Lemma 6.4,
\cite{hessschick98}, p.332)
\begin{eqnarray*}
\lefteqn{\frac{d}{ds}_{|s=0}\left(\frac{1}{\Gamma(s)}\int_0^\ve
k_{P,c}(t)t^{s-1}dt\right)+\int_\ve^\infty k_{P,c}(t)t^{-1}dt}\\
&=&\left(\int_0^\ve
k_{P,c}(t)t^{s-1}dt\right)_{|s=0}+\int_\ve^\infty k_{P,c}(t)t^{-1}dt\\
&=&-2\pi \int_0^c P(y)\:dy\ . 
\end{eqnarray*}
Since $p_l=p_{2u-l}$ we obtain
\begin{prop}\label{mary} 
Let $X=G/K$ with $m(X)=1$, $P=MAN\subset G$ a fundamental parabolic subgroup,
$u=\frac{1}{2}\dim\naaa$, and $X_M=M/K_M$. Then
$$T^{(2)}(X)=(-1)^\frac{n-1}{2}\chi(X^d_M)\frac{\pi Q_X}{\vol(X^d)}\ ,
$$
where
$$ Q_X=\sum_{l=0}^{u-1} (-1)^{l} \int_0^{u-l}p_l(y\cdot\alpha_0)\:dy\ .$$
\end{prop}
In fact, we have proved the proposition for irreducible $X$, but (\ref{hicks})
shows that it holds in the general case, too.

Let us first discuss the case $G=SO(p,q)^0$, $X=X_{p,q}$, $p\le q$ odd. Then $n=pq$ and $X_M=X_{p-1,q-1}$. Furthermore, the polynomials $p_l$ depend only on the complexification of $G$, i.e., on $p+q$.
Thus $Q_{X_{p,q}}=Q_{H^{p+q-1}}$. This proves Proposition \ref{duality} (a).
We emphasize again that $Q_{H^{p+q-1}}$ is a positive rational number \cite{hessschick98}. In fact, Hess-Schick showed that 
$$(-1)^{l} \int_0^{u-l}p_l(y\cdot\alpha_0)\:dy>0\ ,\qquad l=0,\dots, u-1\ .$$

Let $X=SL(3,\R)/SO(3)$. Then $n=5$ and $u=1$. We find that
$$ Q_X=\int_0^1 y^2 dy=\frac{1}{3}\ ,\quad X^d_M=S^2,\ \chi(X^d_M)=2\ .$$
Using for instance (\ref{trotz}) also $\vol(X^d)$ can be easily computed.
This proves Proposition \ref{hesse} and finishes the proof of Theorem \ref{main}.

\end{document}